\documentclass[letter, 10pt, conference]{ieeeconf}      
\IEEEoverridecommandlockouts        

\renewcommand{\baselinestretch}{1}

\usepackage[font=small,skip=0pt]{caption}
\usepackage{etex}
\usepackage[dvipsnames]{xcolor}
\usepackage[vlined,ruled]{algorithm2e}
\usepackage{cite}
\usepackage{booktabs}
\usepackage{array}
\newcolumntype{P}[1]{>{\centering\arraybackslash}p{#1}}

\usepackage{tikz}
\usepackage{pgfplots}
\pgfplotsset{compat=1.6,
        scaled x ticks = false, 
        xticklabel style={/pgf/number format/fixed,/pgf/number format/precision=3},
        } 
\usepackage{standalone}
\usepackage{makecell}

\usepackage{graphicx,color}
\graphicspath{{./figures/}}
\usepackage{epstopdf}

\usepackage{amsmath}
\usepackage{amssymb}
\usepackage{mathrsfs}
\usepackage{dsfont}

\newtheorem{theorem}{Theorem}[section]

\newtheorem{proposition}[theorem]{Proposition}

\newtheorem{assumption}{Assumption}

\def\be{\begin{equation}}
\def\ee{\end{equation}}

\newenvironment{pfof}[1]{\vspace{1ex}\noindent{\itshape Proof of
    #1:}\hspace{0.5em}} {\hfill\QED\vspace{1ex}}

\DeclareMathOperator*{\argmin}{argmin}                                        	

\DeclareSymbolFont{bbold}{U}{bbold}{m}{n}
\DeclareSymbolFontAlphabet{\mathbbold}{bbold}
\newcommand{\vect}[1]{\mathbbold{#1}}
\newcommand{\vones}[1][]{\vect{1}_{#1}}
\newcommand{\vzeros}[1][]{\vect{0}_{#1}}

\newcommand{\real}{\mathbb{R}}

\DeclareMathAlphabet{\mathpzc}{OT1}{pzc}{m}{it}

\newcommand\oprocendsymbol{\hbox{$\square$}}
\newcommand\oprocend{\relax\ifmmode\else\unskip\hfill\fi\oprocendsymbol}

\def\QEDclosed{\mbox{\rule[0pt]{1.3ex}{1.3ex}}} 
\def\QEDopen{{\setlength{\fboxsep}{0pt}\setlength{\fboxrule}{0.2pt}\fbox{\rule[0pt]{0pt}{1.3ex}\rule[0pt]{1.3ex}{0pt}}}}
\def\QED{\QEDopen} 

\renewcommand{\theenumi}{(\roman{enumi})}

\title{Quadratic Performance Analysis of Secondary Frequency Controllers
\thanks{
}}

\author{Bala Kameshwar Poolla, John W. Simpson-Porco, Nima Monshizadeh, and Florian D\"{o}rfler%
    \thanks{%
    J. W. Simpson-Porco is with the Department of Electrical and Computer Engineering, University of Waterloo, ON, Canada. Email: {\tt jwsimpson@uwaterloo.ca}. N. Monshizadeh is with the Engineering and Technology Institute, University of Groningen, The Netherlands. Email: {\tt n.monshizadeh@rug.nl}.
     B. K. Poolla and F. D{\"o}rfler are with the Automatic Control Laboratory, Swiss Federal Institute of Technology (ETH) Z\"urich, Switzerland. Email: {\tt \{bpoolla,dorfler\}@ethz.ch}. The work of B.K. Poolla and F. D\"orfler is supported by ETH Z\"urich funds and the SNF Assistant Professor Energy Grant \#160573. %
}}

\begin{document}
\maketitle
\thispagestyle{empty}
\pagestyle{empty}
 
%
%

\begin{abstract}
This paper investigates the input-output performance of secondary frequency controllers through the control-theoretic notion of $\mathcal{H}_2$ norms. We consider a quadratic objective accounting for the cost of reserve procurement and provide exact analytical formulae for the performance of continuous-time aggregated averaging controllers. Then, we contrast it with distributed averaging controllers-- seeking optimality conditions such as identical marginal costs-- and primal-dual controllers which have gained attention as systematic techniques to design distributed algorithms solving convex optimization problems. Our conclusion is that while the performance of aggregated averaging controllers, such as gather \& broadcast, is independent of the system size and driven predominantly by the control gain, the plain vanilla closed-loop primal-dual controllers scale poorly with size and do not offer any improvement over feed-forward primal-dual controllers. Finally, distributed averaging-based controllers scale sub-linearly with size and are independent of system size in the high-gain limit.

\end{abstract}

\section{Introduction}
\label{Section: Introduction}

Power networks are designed to operate around a nominal frequency (e.g., $50$Hz or $60$Hz) and deviations from this nominal value indicate the global imbalance of power supply and demand. Frequency is, therefore, arguably the most important measurement signal available within the grid, and is taken as a controlled output for controller design over several time-scales to balance supply and demand. To restore the frequency to the nominal operating regime, three hierarchical controls are typically in place. The primary control is a proportional control implemented at sources \cite{JM-JB-JB:11} or loads \cite{CZ-UT-LN-SL:14} and operates over fast time-scales, stabilising the grid to an off-nominal frequency. The secondary control is an integral action which drives the steady-state deviation from the nominal frequency to zero. Such a control is carried out either through centralized approaches such as Automatic Generation Control (AGC) or via decentralized local integral control and does not usually preserve power sharing. The tertiary controls operate on even slower time-scales and aim to minimise the cost of generation, reserves and satisfying operational constraints. This combined optimization problem is referred to as the \emph{optimal frequency regulation problem} (OFRP).

With the expanding integration of distributed generation in the grid, there is a need to revisit the way OFRP is solved as centralised approaches may no longer be efficient in this paradigm. Intermittent renewable sources lead to fast fluctuations in power supply which must be balanced by spinning reserves, which are currently-- expensive, fast-ramping natural gas generators. Frequency control too, needs to be distributed across both the generation and load side to enhance resilience. Controllable power electronic devices offering such regulation need to respond based on minimal information and ideally without detailed knowledge of system parameters or generation/load forecasts. In \cite{NL-CL-ZC-SHL:13} it was shown that the standard primary-secondary control dynamics of a power network can themselves be interpreted as a primal-dual algorithm for solving an OFRP. Other related works include \cite{EM-CZ-SL:17,EM-SL:14,SY-LC:14,XZ-NL-AP:15,TS-CDP-AVDS:17}. The shared characteristics of these primal-dual control strategies are (a) a heavy reliance on power system model information, (b) measurement of power injections/demands or branch-wise power flows, and (c) inter-bus communication of primal or dual variables.

An alternative approach to OFRP based on distributed averaging was proposed in \cite{JWSP-FD-FB:13}, and developed further in the series of papers \cite{JWSP-QS-FD-JMV-JMG-FB:13s,FD-JWSP-FB:16,FD-JWSP-FB:13y}; see also \cite{QS-JMG-JCV:14,MA-DVD-KHJ-HS:13,MA-DVD-HS-KHJ:14,ST-MB-CDP:14,AB-FLL-AD:14,CZ-EM-FD:15,FD-SG:16} for related work on averaging-based secondary frequency controllers. 
Rather than being derived from gradients of a Lagrangian, controllers in this class are consensus-based, `seeking identical marginal costs', and designed such that the desired optimizer is a stable equilibrium point. The shared characteristics among these averaging-based controllers are (a) one or more integral control states, and (b) a consensus term implemented through inter-bus communication driving the source marginal costs towards one another.

In broad strokes then\footnote{For completeness, we note that there are distributed strategies which do not fit cleanly into either of the primal-dual or averaging categories of controllers,  such as the $\lambda$-iteration based controller in \cite{RM-SD-BBC:12}. Here we will focus only on primal-dual and averaging-based controllers.}, we have three large classes of stabilizing controllers which achieve optimal frequency regulation, the gather \& broadcast averaging, the distributed primal-dual method and finally leveraging distributed averaging algorithms. Two natural questions arise: (1) when implemented online as controllers, how well do primal-dual algorithms reject disturbances and (2) how do these disturbance rejection properties compare with the disturbance rejection properties of the averaging-based controllers? Our goal here is to assess and compare the quadratic input-output performance (i.e., $\mathcal{H}_2$) of the OFRP for these three classes of secondary controllers. 

The main contributions of this paper are as follows. We apply and extend some of our previous results in \cite{JWSP-BKP-NM-FD:18, JWSP-BKP-NM-FD:16} to the OFRP for power systems. In Section~\ref{Sec:CentralOFR}, we provide simple and exact analytical formulae for the $\mathcal{H}_2$ performance of closed-loop aggregation-based algorithms subject to load and generation disturbances by considering as controlled output, the cost of the optimization problem. In Section~\ref{Sec:PrimalDualSecondaryPerformance} and Section~\ref{Sec:H2Averaging} we study the $\mathcal{H}_2$ performance of closed-loop primal-dual and distributed averaging-based controllers for the OFRP. Under some simplifying assumptions, we provide exact analytical results showing that in the limit of high gain consensus, the closed-loop system's $\mathcal{H}_2$ performance under distributed averaging-based controllers becomes completely independent of system dimensionality, converging to the performance of the broadcast controllers. In Section \ref{Section: Discussion} we provide a detailed discussion and comparison of primal-dual, distributed averaging-based, and broadcast controllers. Then, we take a numerical approach and examine system performance for more general objective functions. The results confirm the conclusions that distributed averaging-based, broadcast controllers offer superior disturbance rejection than primal-dual controllers.

\section{Power Network Modeling and $\mathcal{H}_2$ Norms}
\label{sec:modeling}

\subsection{Power Network Modeling}
We model a power network as a weighted graph $(\mathcal{V},\mathcal{E})$ where $\mathcal{V} = \{1,\ldots,n\}$ is the set of nodes (buses) and $\mathcal{E} \subset \mathcal{V}\times\mathcal{V}$ is the set of edges (branches) with associated edge weights $b_{ij} > 0$ for $\{i,j\} \in \mathcal{E}$. For each bus $i \in \mathcal{V}$ we associate state variables $(\theta_i,\omega_i)$ corresponding to the voltage phase angle and the frequency deviation from nominal. Under the linear DC Power Flow approximation, the system evolves according to the swing dynamics
\begin{equation}
\label{Eq:Swing}
\begin{aligned}
\begin{bmatrix}
\dot{\theta}\\
M\dot{\omega}
\end{bmatrix}
&=
\begin{bmatrix}
\vzeros[] & I_n \\ 
-L & -D 
\end{bmatrix}\begin{bmatrix}\theta \\ \omega \end{bmatrix} + 
\begin{bmatrix}
\vzeros[]\\
P^\star + p
\end{bmatrix},
\end{aligned}
\end{equation}
where $M={\rm diag} (m_i)$, $D={\rm diag} (d_i)$, $m_i > 0$ represents inertia or inverter filter time constant, $d_i > 0$ models damping and/or droop control, $P^\star$ is vector of constant nominal active power injections, $L$ is the network Laplacian, $p$ is the vector of control inputs corresponding to additional power generation from reserves, and $\vzeros[]$ is a matrix/vector of zeros of appropriate dimension. Throughout, we ignore reactive power and voltage dynamics; these assumptions are standard in secondary frequency control studies.

When $p = \vzeros[n]$, the dynamics \eqref{Eq:Swing} converge from every initial condition to a common steady-state frequency $\omega \rightarrow \omega_{\rm{ss}}\vones[n]$ which can be easily calculated to be $\omega_{\rm{ss}} = ({\sum_{i=1}^n P_i^\star})/({\sum_{i=1}^n d_i})$. When $\omega_{\rm{ss}} \neq 0$, this represents a static deviation from nominal, which we will eliminate by appropriately selecting the reserve secondary power inputs $p$. Moreover, since the variables $P_i^\star$ will tend to vary, feedback control should be used to provide real-time frequency regulation. To determine the steady-state values for $p_i$, an OFRP can be formulated as%
\begin{subequations}
\label{Eq:DAPIOptimizationScalar}
\begin{alignat}{2}
\label{Eq:DAPIOptimizationScalar1}
& \underset{{p \in \real^n}}{\text{minimize}} &\qquad & \mathlarger{\sum_{i=1}^n} \frac{1}{2}k_i p_i^2 \\\label{Eq:DAPIOptimizationScalar2}
& \text{subject to} & & 0 = \vones[n]^{\top}(P^\star + p)\,,
\end{alignat}
\end{subequations}
where we seek to minimize the total cost\footnote{A linear term could also of course be added to the cost, but we omit it here for simplicity. We assume that any inequality constraints on $p$ are non-binding, and subsequently drop them from the problem \eqref{Eq:DAPIOptimizationScalar}, which is then similar to the \emph{classic} economic dispatch (see \cite{ARB:09}).} \eqref{Eq:DAPIOptimizationScalar1} of reserve generation $p_i \in \real$, for some cost coefficients $k_i > 0$. The minimization is subject to network-wide balancing of power injections \eqref{Eq:DAPIOptimizationScalar2}.
One may deduce from \eqref{Eq:Swing} that the constraint \eqref{Eq:DAPIOptimizationScalar2} also ensures that $\omega = \vzeros[n]$ in steady-state, i.e., the frequency returns to its nominal value. The vector of optimal power inputs can be computed as $p_{\rm opt}=-(\vones[n]^{\top}P^\star/\vones[n]^{\top}K^{-1}\vones[n]) K^{-1} \vones[n]$, where $K={\rm diag}(k_i)$.

\subsection{System Performance in the $\mathcal{H}_2$ Norm}
We very briefly review $\mathcal{H}_2$ norms of a linear system, which is used to evaluate the performance of the secondary frequency controllers. 
Consider the MIMO continuous-time LTI system, $G=(A, B, C)$ with a state-space description
\begin{equation}\label{Eq:GeneralMIMO}
\begin{aligned}
\dot{x} &= Ax + B\eta\\
y &= Cx\,,
\end{aligned}
\end{equation}
where $\eta$, $y$ are the input, output respectively, and $A$ is Hurwitz. If \eqref{Eq:GeneralMIMO} is input-output stable, its $\mathcal{H}_2$ norm, $\|G\|_{\mathcal{H}_2}$ measures the steady-state output variation under stochastic disturbances, i.e., 
\[ \|G\|_{\mathcal{H}_2}^2 ={\lim_{t \to \infty}} \mathbb{E}[y(t)^\top y(t)]\,,\]
when each component of $\eta(t)$ is stochastic white noise with unit covariance (i.e., $\mathbb{E}[\eta(t)\eta(\tau)^\top]=\delta(t-\tau)I$). A convenient formula for the $\mathcal{H}_2$ norm is \cite[Chapter 6]{GED-FP:00}
\begin{equation}\label{Eq:H2Obsv}
\|G\|_{\mathcal{H}_2}^2 = \mathrm{Tr}(B^\top XB)\,,
\end{equation}
where $X = X^\top> 0$ is the observability {Gramian} satisfying
\begin{equation}\label{Eq:ObsvGram}
XA + A^\top X + C^\top C = \vzeros[].
\end{equation}
If the pair $(C,A)$ is observable, then \eqref{Eq:ObsvGram} is solvable for the unique, positive-definite observability {Gramian}. 

In what follows we will derive results on the $\mathcal{H}_2$ performance for different algorithms under disturbances by explicitly solving \eqref{Eq:ObsvGram} for particular cases. Recent applications of the $\mathcal{H}_2$ norm to power system performance may be found in \cite{ET-BB-DFG:15,ET-MA-JWSP-HS:15d,BKP-SB-FD:16,ERAW-YJ-CZ-EM-CDP-FD:17}.

\section{Broadcast Algorithm Performance for Optimal Frequency Regulation}
\label{Sec:CentralOFR}

In this section we examine the $\mathcal{H}_2$ performance of closed-loop centralized algorithms for the OFRP \eqref{Eq:DAPIOptimizationScalar}. In particular, we examine the performance of the centralized frequency averaging algorithm. These centralized results will serve as useful benchmarks for later comparison to results derived for distributed algorithms.

Beginning from the OFRP \eqref{Eq:DAPIOptimizationScalar} and partially following \cite{CZ-UT-LN-SL:14}, the Lagrangian of \eqref{Eq:DAPIOptimizationScalar} is given by
$$
L^g (p,\mu) = \frac{1}{2}p^{\top}Kp + \mu\vones[n]^{\top}(P^\star+p)\,,
$$
where $\mu \in \real$ is a multiplier. On computing the $\argmin_p L^g(p,\mu)$, one finds that $p_{\rm opt} =-\mu K^{-1}\vones[n]$ is the unique minimizer.
Next, we assume the existence of a central aggregator \cite{QS-JMG-JCV:14,FD-JWSP-FB:16,FD-SG:16} which collects frequency measurements from each bus in the network. The aggregator averages the frequency errors, integrates the result, and broadcasts control signals back to each bus via%
\begin{subequations}
\label{Eq:CAPIComponents}
\begin{align}
\tau_\mu \dot{\mu}&=\sum_{j=1}^{n}\nolimits r_j \omega_j\,,\\
p &= -\mu K^{-1}\vones[n]\,
\end{align}
\end{subequations}
where $\tau_\mu>0$ is a scalar gain and $\{r_i\}_{i=1}^n$ are a set of convex coefficients. Note that the dynamics in \eqref{Eq:CAPIComponents} encapsulate Automatic Generation Control (ACG), Centralised Averaging PI (CAPI), and Gather \& Broadcast schemes. Now consider the interconnection of the central averaging controller \eqref{Eq:CAPIComponents} with the power system dynamics \eqref{Eq:Swing}. One may verify without much effort that $(\theta_{\rm opt},\omega_{\rm opt},\mu_{\rm opt}) = (\theta_{\rm opt},\vzeros[n],\vones[n]^{\top}P^\star/(\vones[n]^{\top}K^{-1}\vones[n]))$ is the unique\footnote{Up to a uniform translation of all phase angles.} stable equilibrium point of the interconnection and $\mu_{\rm opt}$ is optimal for \eqref{Eq:DAPIOptimizationScalar}. We therefore obtain frequency regulation. To study the input-output performance, we translate this equilibrium to the origin and assume that the power injection vector $P^\star$ is corrupted by an additive noise (modeling fluctuating generation/
load, noise, or other uncertainties) $B_1\eta$, for $B_1 \in \real^{n \times n}$. As the performance output, we select $\|y(t)\|_2^2 = p^{\top}Kp$. Under this scenario, the closed-loop input-output dynamics, with the transformation $\varphi=E^\top \theta$ in vector form is
\begin{equation}
\label{Eq:CAPI}
\begin{aligned}
\begin{bmatrix}
\dot{\varphi}\\
M\dot{\omega}\\
\tau_\mu \dot{\mu}\\
\end{bmatrix}
&=
\begin{bmatrix}
\vzeros[] & E^{\top} & \vzeros[n-1]  \\ 
-E & -D & -K^{-1}\vones[n] \\ 
\vzeros[n-1]^{\top} & \vones[n]^{\top}R  & 0 \\
\end{bmatrix}\begin{bmatrix}\varphi \\ \omega \\ \mu\end{bmatrix} + 
\begin{bmatrix}
\vzeros[]\\
B_1\\
\vzeros[]
\end{bmatrix}\eta,\,\\
y &= -(K^{-\frac{1}{2}}\vones[n])\mu\,,
\end{aligned}
\end{equation}
where $E$ is the incidence matrix, $R=\mathrm{diag}(r_1,\ldots,r_n)$.
\smallskip
\begin{assumption}[Uniform Parameters]\label{Assumption:Parameters} Throughout the paper we assume uniform parameters, i.e., $M = m I_n$, $D = d I_n$, $B_1=b I_n$, and $R= \frac{1}{n} I_n$. While the assumption is restrictive, it allows us to derive some benchmark results.
\end{assumption}
\smallskip
\begin{assumption}[Communication Topology]\label{Assumption:Communication} We assume that the communication through which the controller is implemented has the same structure as the underlying electrical network modulo a constant factor $\gamma>0$, i.e., $L_{\rm c}=E_{\rm c} E_{\rm c}^\top=\gamma E E^\top=\gamma L$. Furthermore, we assume the the underlying network is acyclic, strictly for simplicity. This assumption can be relaxed without affecting the conclusions.
\end{assumption}
\smallskip

\begin{assumption}[Identical Costs]\label{Assumption:Costs} We assume that the cost of reserve procurement is identical, i.e., $K=k I_n$.
\end{assumption}
\smallskip

The following result characterizes the $\mathcal{H}_2$ norm of the system for a simplifying choice of parameters.
\smallskip
\begin{theorem}[Performance of Broadcast Algorithm]
\label{Thm:CAPIPerformance} We consider the centralized averaging OFRP dynamics \eqref{Eq:CAPI} with disturbances $\eta$ and performance outputs $y$.  Under Assumptions (1)--(3) and for some $\tau_\mu, b >0$, we have that the squared $\mathcal{H}_2$ norm of \eqref{Eq:CAPI} is
\begin{equation}
\begin{aligned}\label{Eq:FrequencyPMH2CostOutput1}
\|G\|_{\mathcal{H}_2}^2 &= \frac{b^2}{2 \tau_\mu}\frac{1}{d} \,.
\end{aligned}
\end{equation}
\end{theorem}
\smallskip
Please refer to the Appendix for the proof.
\smallskip

Note that the broadcast averaging algorithm does not scale with the system size, further for sufficiently large networks and moderate damping coefficients, the performance improvement in terms of disturbance rejection over the plain vanilla gradient ascent algorithm in \cite{JWSP-BKP-NM-FD:16} is significant.

\section{Distributed Algorithm Performance for Optimal Frequency Regulation}
\label{Sec:ApplicationToOFR}

In this section we assess the performance of distributed primal-dual and averaging-based frequency controllers for power systems. 

\subsection{$\mathcal{H}_2$ Performance of Primal-Dual Frequency Controllers}
\label{Sec:PrimalDualSecondaryPerformance}
We recall from the Lagrangian of \eqref{Eq:DAPIOptimizationScalar} that the dual function $\Phi(\mu) = \inf_{p \in \real^n}L^g (p,\mu)$, and the OFRP becomes
\begin{equation}\label{Eq:AGC}
\underset{\mu \in \real}{\text{maximize}} \quad \Phi(\mu) = \sum_{i=1}^n \mu P_i^\star - \frac{1}{2k_i}\mu^2\,,
\end{equation}
where we seek to maximize $\Phi(\mu)$ over $\mu \in \real$. On introducing local variables $\mu_i \in \real$ for each bus and a consensus constraint, the problem \eqref{Eq:AGC} is again equivalent to
\begin{equation}
\label{Eq:ZhaoProblem}
\begin{aligned}
& \underset{\mu \in \real^n}{\text{maximize}}
& & \sum_{i=1}^n\nolimits \mu_iP_i^\star - \frac{1}{2k_i}\mu_i^2 \\
& \text{subject to}
&&  0= \mu_i - \mu_j \,, \qquad \{i,j\}\in\mathcal{E}_{\rm c}\,,
\end{aligned}
\end{equation}
where $\mathcal{E}_{\rm c}$ is the edge set of a connected, undirected, and acyclic communication graph $(\mathcal{V},\mathcal{E}_{\rm c})$ between the buses. The additional constraints $\mu_i - \mu_j = 0$ force the local variables $\mu_i$ to agree at optimality. By letting $E_{\rm c} \in \real^{n \times |\mathcal{E}_{\rm c}|}$ denote the incidence matrix \cite{CAD-ESK:69} of the communication graph, the dual OFRP \eqref{Eq:ZhaoProblem} is written in vector notation as
\begin{equation}\label{Eq:ZhaoProblemV}
\begin{aligned}
& \underset{\mu \in \real^n}{\text{minimize}}
& & \frac{1}{2}\mu^{\top}K^{-1}\mu - (P^\star)^{\top}\mu \\
& \text{subject to}
&  &\vzeros[|\mathcal{E}_{\rm c}|]= E_{\rm c}^{\top}\mu \,,
\end{aligned}
\end{equation}
where now $\mu = (\mu_1,\ldots,\mu_n)$. The problem \eqref{Eq:ZhaoProblemV} is a  linearly constrained, strictly convex quadratic program. The corresponding Lagrangian is $L^g(\mu,\nu) = \frac{1}{2}\mu^{\top}K^{-1}\mu - (P^\star)^{\top}\mu + \nu^{\top}E_{\rm c}^{\top}\mu$, where $\nu \in \real^{|\mathcal{E}_{\rm c}|}$ is a vector of multipliers. The primal-dual algorithm \cite{JWSP-BKP-NM-FD:18} with control $p$ thus becomes

\begin{equation}\label{Eq:PrimalDualSecondaryDynamics}
\begin{aligned}
\mathcal{T}_{\mu}\dot{\mu} &= -K^{-1}\mu + P^\star - E_{\rm c}\nu\,,\\
\mathcal{T}_{\nu}\dot{\nu} &= E_{\rm c}^{\top}\mu,\quad p= -K^{-1}\mu\,,
\end{aligned}
\end{equation}
where $\mathcal{T}_{\mu}$ and $\mathcal{T}_\nu$ are positive diagonal matrices of controller gains. We note that the algorithm \eqref{Eq:PrimalDualSecondaryDynamics} does not make use of any real-time state information from \eqref{Eq:Swing}, but instead requires knowledge of the nominal power injections $P^\star$ to determine the optimal feed-forward set-point for $p$. 

A common variation (see, for example, \cite{TS-CDP-AVDS:17}) on the primal-dual frequency controller \eqref{Eq:PrimalDualSecondaryDynamics} is to add frequency deviation feedback to the state equation for $\mu$ as 
\begin{equation}\label{Eq:TauEquationModified}
\mathcal{T}_{\mu}\dot{\mu} = -K^{-1}\mu + P^\star - E_{\rm c}\nu + \alpha K^{-1}\omega\,,
\end{equation}
where $\alpha > 0$ is a gain (typically equal to one). In vector form, the closed-loop dynamics now read as
\begin{equation}
\label{Eq:PrimalDual_Feedback_IO}
\begin{aligned}
\begin{bmatrix}
\dot{\varphi}\\
M\dot{\omega}\\
\mathcal{T}_\mu \dot{\mu}\\
\mathcal{T}_\nu \dot{\nu}
\end{bmatrix}
&\!=\!
\begin{bmatrix}
\vzeros[] & E^\top & \vzeros[] & \vzeros[] \\ 
-E & -D & -K^{-1} & \vzeros[] \\ 
\vzeros[] & \alpha K^{-1} & -K^{-1} & -E_{\rm c}\\
\vzeros[] & \vzeros[] & E_{\rm c}^{\top} & \vzeros[]
\end{bmatrix}\begin{bmatrix}\varphi \\ \omega \\ \mu \\ \nu\end{bmatrix}\! +\! 
\begin{bmatrix}
\vzeros[]\\
B_1\\
B_1\\
\vzeros[]
\end{bmatrix}\eta\\
y &= {-}K^{-\frac{1}{2}}\mu,
\end{aligned}
\end{equation}
where the output is same as in \eqref{Eq:CAPI}.
We denote the associated input-output map of this modified system \eqref{Eq:PrimalDual_Feedback_IO} by $G(\alpha)$. When we set $\alpha = 0$, the system \eqref{Eq:PrimalDual_Feedback_IO} reduces to the cascade of the primal-dual algorithm with the power system dynamics \eqref{Eq:Swing}, and for $\mathcal{T}_\mu=\tau I_n$, $\mathcal{T}_\nu=\tau I_n$, the $\mathcal{H}_2$ norm $\|G(0)\|_{\mathcal{H}_2}^2=(b^2/2\tau) n$ \cite[Theorem 4.1]{JWSP-BKP-NM-FD:16}.  While naively one may expect that adding real-time frequency information to the algorithm should improve performance, observe that a non-zero $\alpha$ term introduces additional skew-symmetric structure in the system matrix of \eqref{Eq:PrimalDual_Feedback_IO} {(for $\alpha=0$ \eqref{Eq:PrimalDual_Feedback_IO} is a simple cascade system)}, suggesting additional oscillations and therefore, worse transient performance. 
\smallskip

\begin{theorem}[Primal-Dual Performance]\label{Thm:PrimalDualFeedbackPerformance}
We now consider the closed-loop primal-dual dynamics \eqref{Eq:PrimalDual_Feedback_IO} with disturbance inputs $\eta$ and performance outputs $y$, and with additional frequency feedback. We assume that Assumptions (1)--(3) hold and that $\mathcal{T}_\nu$, and $\mathcal{T}_\mu$ are multiples of the identity matrix. Then for any $\alpha \geq 0$, the $\mathcal{H}_2$ norm of the system \eqref{Eq:PrimalDual_Feedback_IO} satisfies 
$$\|G(\alpha)\|_{\mathcal{H}_2}^2 \leq \dfrac{b^2}{2 \tau} n + \frac{b^2\alpha n}{2m}.$$
\end{theorem}
\smallskip
Please refer to the Appendix for the proof.
\smallskip

Theorem \ref{Thm:PrimalDualFeedbackPerformance} indicates that in the special case considered \textemdash{} the  additional frequency feedback considered in \eqref{Eq:TauEquationModified} does not improve controller performance, and in fact worsens it compared to the original feed-forward design \eqref{Eq:PrimalDualSecondaryDynamics}. We refer the reader to Table~\ref{tab:case} for a numerical analysis.

\subsection{$\mathcal{H}_2$ Performance of Distributed Averaging-Based Frequency Controllers}
\label{Sec:H2Averaging}

An alternative approach for optimal frequency regulation of power systems was proposed in \cite{JWSP-FD-FB:13}, based on distributed averaging algorithms. To begin, we return to the OFRP \eqref{Eq:DAPIOptimizationScalar}, and we note that the constraint \eqref{Eq:DAPIOptimizationScalar2} is satisfied if and only if there exists a vector $v \in \real^{n}$ such that $\vones[n]^{\top}v = 0$ and $P^\star + p - v = \vzeros[n]$.  Moreover, any such $v$ is in the image of the network Laplacian matrix $L$, and hence can be written (with some notational foresight) as $v = L\theta$ for $\theta \in \real^n$. The OFRP \eqref{Eq:DAPIOptimizationScalar} is therefore equivalent to 

\begin{subequations}\label{Eq:DAPIOptimization}
\begin{align}
\label{Eq:DAPIOptimization1}
& \underset{\theta \in \mathbb{R}^n,\, p \in \real^n}{\text{minimize}}
& & \sum_{i=1}^n\nolimits \frac{1}{2}k_ip_i^2 \\
\label{Eq:DAPIOptimization2}
& \text{subject to}
&  \vzeros[n] &= P^\star -L\theta + p\,.
\end{align}
\end{subequations}
Comparing the dynamic equations \eqref{Eq:Swing} and the hard constraint \eqref{Eq:DAPIOptimization2}, we see that in steady-state $\omega = \vzeros[n]$, and therefore frequency regulation is achieved. 
The optimal points of \eqref{Eq:DAPIOptimization} can in fact be exactly determined as follows, the proof of which follows from applying the KKT conditions.

\smallskip

\begin{proposition}[Optimal Secondary Inputs and States]
\label{Prop:OFRPKKT}
The unique\footnote{Up to a uniform rotation of all phase angles $\theta$.} primal optimizer $(p_{\rm opt},\theta_{\rm opt})$ of the OFRP \eqref{Eq:DAPIOptimization} is given by
\begin{equation}\label{Eq:OFRPOptimizers}
p_{\rm opt} = -\frac{\vones[n]^{\top}P^\star}{\vones[n]^{\top}K^{-1}\vones[n]}K^{-1}\vones[n], \, \theta_{\rm opt} = L^{\dagger}(P^\star+p_{\rm opt}),\,
\end{equation}
\end{proposition}
where $L^\dagger$ is the pseudo-inverse of the network Laplacian $L$.

In contrast to primal-dual algorithms which derive from gradients of the desired cost function, averaging-based strategies are carefully designed to converge to the optimal points \eqref{Eq:OFRPOptimizers} such as identical marginal costs. The secondary control input $p=K^{-1}\lambda$ in \eqref{Eq:Swing} is designed component-wise as
\begin{equation}\label{Eq:DAPISecondaryComponents}
\tau_i\dot{\lambda}_i = -\omega_i - \sum_{j=1}^n \nolimits a_{ij}(\lambda_i-\lambda_j),
\end{equation}
where $\tau_i$ is a controller gain, $k_i$ is the cost coefficient from the optimization problem \eqref{Eq:DAPIOptimization}, and the coefficients $a_{ij} = a_{ji} \geq 0$ are the elements of a symmetric adjacency matrix $\mathcal{A}_{\rm c} = \mathcal{A}_{\rm c}^{\top} \in \real^{n\times n}$ describing a weighted, undirected and connected communication network $(\mathcal{V},\mathcal{E}_{\rm c})$ among the buses. The first term in \eqref{Eq:DAPISecondaryComponents} integrates the frequency error, while the second term uses inter-bus communication to drive the unconstrained marginal costs $k_ip_i$ for each bus towards one another. In vector notation, the controller \eqref{Eq:DAPISecondaryComponents} reads as

\begin{equation}\label{Eq:DAPISecondary}
\mathcal{T}K\dot{p} = -\omega - L_{\rm c}K\lambda,
\end{equation}
where $\mathcal{T}=\tau I_n$ is a diagonal matrix of controller gains and $L_{\rm c} = \mathrm{diag}(\mathcal{A}_{\rm c}\vones[n])-\mathcal{A}_{\rm c}=\gamma L$ is a symmetric Laplacian matrix for $\gamma>0$. The following result indicates convergence of these dynamics to the global minimizer of \eqref{Eq:DAPIOptimization}; see \cite{FD-JWSP-FB:13y} for a proof using first-order inverter dynamics.

\smallskip

As we did in Section \ref{Sec:PrimalDualSecondaryPerformance} for the primal-dual frequency controller, we now study the input-output performance of the distributed averaging controller when the nominal power injections $P^\star$ are subject to disturbances. We again shift to error coordinates, and as a performance output we take the controller state $p$ with weighting matrix $K$ from the optimization problem \eqref{Eq:DAPIOptimization} by selecting $y= K^{\frac{1}{2}}p$, such that $\|y(t)\|_{2}^2 = p^{\top}Kp$. 
The corresponding input-output dynamics of \eqref{Eq:Swing}, \eqref{Eq:DAPISecondary} are
\begin{equation}\label{Eq:DAPI_IO}
\begin{aligned}
\begin{bmatrix}
\dot{\theta}\\
M\dot{\omega}\\
\mathcal{T}K\dot{p}
\end{bmatrix}
&=
\begin{bmatrix}
\vzeros[] & I_n & \vzeros[] \\ 
-L & -D & I_n \\ 
\vzeros[] & -I_n & -L_{\rm c}K
\end{bmatrix}\begin{bmatrix}\theta \\ \omega \\ p\end{bmatrix} + 
\begin{bmatrix}
\vzeros[]\\
B_1\\
\vzeros[]
\end{bmatrix}\eta\\
y &= K^{\frac{1}{2}}p,
\end{aligned}
\end{equation}
where the output is same as in \eqref{Eq:CAPI}. In general, computing an insightful analytic formula for the $\mathcal{H}_2$ norm of the system \eqref{Eq:DAPI_IO} is difficult. Under Assumptions (1)--(3), however, the norm can be computed in closed-form. 

\smallskip

\begin{theorem}[Distributed Averaging Performance]\label{Thm:DAPI_Performance}
We consider the closed-loop power system \eqref{Eq:DAPI_IO} under the distributed averaging-based controller \eqref{Eq:DAPISecondary}, with performance output $\|y(t)\|^{2}_2 = p^{\top}Kp$ and disturbance input $\eta \in \real^n$. Then the squared $\mathcal{H}_2$ norm of the system \eqref{Eq:DAPI_IO} is 
\begin{equation}\label{Eq:DAPIGeneralH2Norm}
\|G\|_{\mathcal{H}_2}^2 = \frac{b^2}{2\tau}\frac{1}{d}\sum_{i=1}^n \frac{1}{z_2\lambda_i^2 + z_1\lambda_i + 1}\,,
\end{equation}
where $z_2 = mk\gamma^2/\tau$, $z_1 = m\gamma/(d\tau) + kd\gamma + k\tau$, and $\lambda_i$ is the $i$th eigenvalue of the grid Laplacian matrix $L$. Moreover, the following special cases hold:
\begin{enumerate}
\item \textbf{Overdamped Limit:} It holds that
\begin{equation}\label{Eq:DAPI_H2_Inverters}
\lim_{m \rightarrow 0}\|G\|_{\mathcal{H}_2}^2 = \frac{b^2}{2\tau}\frac{1}{d}\sum_{i=1}^n \frac{1}{1+(k\tau+kd\gamma)\lambda_i}\,;
\end{equation}
\item \textbf{High Averaging Gain:} It holds that 
\begin{equation}\label{Eq:DAPI_H2_HighGain}
\lim_{\gamma \rightarrow \infty} \|G\|_{\mathcal{H}_2}^2 = \frac{b^2}{2\tau}\frac{1}{d}\,.
\end{equation}
\end{enumerate}
\end{theorem}
\smallskip
Please refer to the Appendix for the proof.
\smallskip

While the general formula \eqref{Eq:DAPIGeneralH2Norm} is somewhat opaque, the overdamped limit \eqref{Eq:DAPI_H2_Inverters} is considerably simpler, and applies when generators have small inertias or when measurement and actuation delays at inverters are negligible \cite{JS-DZ-RO-AMS-TS-JR:16}. 
In both the general case \eqref{Eq:DAPIGeneralH2Norm} and the overdamped limit \eqref{Eq:DAPI_H2_Inverters}, the squared $\mathcal{H}_2$ norm is given by a sum of network modes, with the terms in the sum scaling inversely with the eigenvalues of the grid Laplacian. 
Since in realistic networks the Laplacian eigenvalues $\{0, \lambda_2,\lambda_3,\lambda_4,\ldots\}$ are a rapidly increasing sequence, this implies that the squared $\mathcal{H}_2$ norm \eqref{Eq:DAPIGeneralH2Norm} scales sub-linearly with system size. 

Of particular interest is the fact that the terms in the sum scale inversely with $\gamma$, the gain on the averaging term in the controller. The special case (ii) indicates that the $\mathcal{H}_2$ norm can be suppressed by increasing $\gamma$ sufficiently. Note in particular that in the high-gain limit \eqref{Eq:DAPI_H2_HighGain}, the $\mathcal{H}_2$ norm is \emph{independent of system size} and identical to the broadcast performance \eqref{Eq:FrequencyPMH2CostOutput1}. This appears to be a fundamental difference between the averaging-based controllers and the open-loop primal-dual controller \eqref{Eq:PrimalDualSecondaryDynamics}. Intuitively, the high averaging gain limit forces the agents to quickly agree on their unconstrained marginal costs $k_ip_i$, which immediately minimizes the cost \eqref{Eq:DAPIOptimizationScalar1} in the absence of the equality constraint \eqref{Eq:DAPIOptimizationScalar2}. The remaining slow dynamics for $p$ are one-dimensional, driving the system towards feasibility of the equality constraint \eqref{Eq:DAPIOptimizationScalar2}, and hence to frequency regulation.

\section{Discussion and Simulations}
\label{Section: Discussion}

We now make some general comments comparing the three classes of controllers, delineating their relative advantages and disadvantages. 

The primal-dual based controllers are quite flexible, in that they can be modified to handle hard inequality constraints and strictly convex objective and lead to a systematic design procedure. However, they typically require model information such as inertia and damping coefficients, line susceptances, or in the case of \eqref{Eq:PrimalDualSecondaryDynamics}, disturbance measurements. The resulting algorithms also tend to assign dynamic states to edges of the network, requiring additional processors for distributed implementation. The primal-dual controller \eqref{Eq:PrimalDualSecondaryDynamics} is an online method to compute the optimal feed-forward input $p$ given measurements of the disturbances affecting $P^\star$, suggesting that it is potentially non-robust to modeling errors and parametric uncertainties. It remains unclear to what extent the modification \eqref{Eq:TauEquationModified} remedies this, or how communication delays will influence algorithm stability. 

The broadcast-based controllers on the other hand are comparatively simpler to understand. They rely on a single system-wide integrator to average the frequency errors and broadcast the control signals back to each bus. They do not need model information for implementation and are true feedback controllers measuring only the local frequency and produces a corresponding control input, requiring no disturbance measurement. Though easier to implement and tune, the presence of only one controller renders them less robust to communication malfunctions.

Distributed averaging strategies such as \eqref{Eq:DAPISecondaryComponents} are a concatenation of an integral control and an averaging of the marginal costs. The implementation is comparatively simpler than the primal-dual controllers, requiring only one control state to be assigned to each bus. Similar to broadcast controllers, they too do not need model information. Moreover, the consensus-type algorithms on which these controllers are founded are known to be quite robust to communication malfunctions such as asynchronism, delays, and packet losses \cite{VDB-JMH-AO-JNT:05}. The main disadvantage of distributed averaging controllers is that it is yet unknown how to accommodate hard constraints. Moreover, convergence of the algorithm for non-quadratic cost functions (see \cite[Equation 11]{FD-SG:16}) in the OFRP remains an open research problem. 

Finally, primal-dual controllers tend to suffer from scalability issues (see \cite{JWSP-BKP-NM-FD:16} for augmented Lagrangian implementation to counteract this) and feedback controllers based on averaging approaches perform much better in large networks.

\begin{table}[h!]
\begin{center}
\caption{The squared $\mathcal{H}_2$ norm for the primal-dual, broadcast, and distributed controllers with a varying penalty on $\omega$.\\}
\label{tab:case}
\bgroup
\def\arraystretch{1.35}
\setlength{\tabcolsep}{1pt}
\begin{tabular}{P{1.5cm} P{1.5cm} P{1.5cm} P{1.5cm} P{1.5cm}}
\hline\noalign{\smallskip}
\multicolumn{1}{c}{$\omega$ penalty} &  {Primal-Dual}   & Primal-Dual & Distributed & Broadcast\\
\multicolumn{1}{c}{{$\sqrt \pi$ }} &  {$\alpha=0$}   & $\alpha=5$ & $\gamma=5$ &\\
\noalign{\smallskip}\hline\noalign{\smallskip}
$0.0$ & $0.417$ & $0.569$ & $0.088$ & $0.083$ \\
$0.3$ & $0.639$ & $0.791$ & $0.311$ & $0.308$ \\
$0.6$ & $1.307$ & $1.458$ & $0.981$ & $0.983$ \\
$0.9$ & $2.421$ & $2.569$ & $2.095$ & $2.108$ \\
$1.2$ & $3.980$ & $4.125$ & $3.656$ & $3.683$ \\
$1.5$ & $5.984$ & $6.125$ & $5.663$ & $5.708$ \\
\noalign{\smallskip}\hline\noalign{\smallskip}
\end{tabular}
\egroup
\end{center}
\end{table}
\vspace{-1em}
To numerically examine the various controllers discussed in the paper, we consider an acyclic network of $5$ nodes and with identical parameters-- $m=1$, $d=1$, $\tau_\mu=\tau_\nu=6$, $k=4$, and $b=1$. Though we do not explicitly consider the system frequency in the performance output $y$, we report the numerical observations in Table~\ref{tab:case}. The performance $\|G\|_{\mathcal{H}_2}^2$ of the three controllers, with a modified objective $y={\rm diag}(K^{\frac{1}{2}}p, \Pi^{\frac{1}{2}}\omega)$ for $\Pi=\pi I_n$ is tabulated as a function of the frequency penalty $\sqrt\pi$. Note that as we increase the feedback gain $\alpha$, the primal-dual performance degrades. Further, the broadcast and the distributed algorithm ($\gamma=5$) have almost identical disturbance rejection.

\begin{figure}[htbp]
\centering
\input{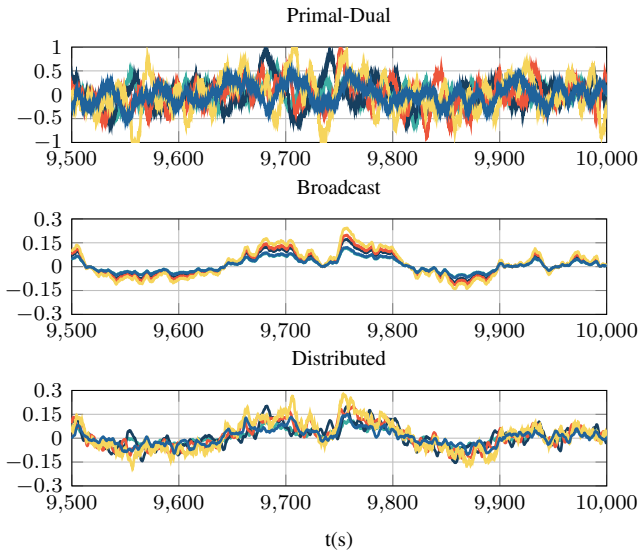}
\caption{The time series (steady-state) of the performance outputs for the three implementations (note the scale difference)-- primal-dual ($\alpha=1$), broadcast, and distributed ($\gamma=1$) controllers for white noise inputs and non-identical $m$, $k$, $d$ parameters.}
\label{Fig:timedom}
\end{figure}

In Figure~\ref{Fig:timedom}, we consider non-uniform parameters for $M$, $D$, $K$. Note that the system performance (measured as the steady-state variance of the output $y$ for white noise input $\eta$) is in consonance with the theoretical results derived previously for uniform parameters.

\vspace{-0.05em}

\section{Conclusions}
\label{Section: Conclusions}

We have quantified the input-output performance of secondary frequency controllers of power systems and compared the $\mathcal{H}_2$ performance of continuous-time primal-dual, distributed averaging, and broadcast methods providing explicit formulae. Our findings indicate that the performance of plain vanilla primal-dual controllers with frequency feedback scales poorly with system size and does not improve performance over the feed-forward open-loop implementation. Further, for aggregated averaging (broadcast) algorithms, the performance is predominantly determined by the control gains and independent of the system size. Standard distributed averaging algorithms, seeking optimality conditions such as identical marginal costs suffer from sub-linear dependence on system size. This can be improved by high control gain, which retrieves the broadcast performance.

An important extension of this work would be to consider objective functions other than control effort or cost of reserve generation-- for example investigate the system performance for objectives involving frequency excursions or a combination thereof theoretically (we provide numerical results). It is also currently unclear how the present results change when including measurement noise in the frequency measurements or when considering higher-order generator dynamics such as governors, frequency-dependent loads, and static loads.

\section*{Acknowledgments}
The authors thank Tjerk Stegink for his input on primal-dual frequency control stability and design.


\renewcommand{\baselinestretch}{1}
\bibliographystyle{IEEEtran}
\bibliography{final}

\begin{thebibliography}{10}
\providecommand{\url}[1]{#1}
\csname url@samestyle\endcsname
\providecommand{\newblock}{\relax}
\providecommand{\bibinfo}[2]{#2}
\providecommand{\BIBentrySTDinterwordspacing}{\spaceskip=0pt\relax}
\providecommand{\BIBentryALTinterwordstretchfactor}{4}
\providecommand{\BIBentryALTinterwordspacing}{\spaceskip=\fontdimen2\font plus
\BIBentryALTinterwordstretchfactor\fontdimen3\font minus
  \fontdimen4\font\relax}
\providecommand{\BIBforeignlanguage}[2]{{%
\expandafter\ifx\csname l@#1\endcsname\relax
\typeout{** WARNING: IEEEtran.bst: No hyphenation pattern has been}%
\typeout{** loaded for the language `#1'. Using the pattern for}%
\typeout{** the default language instead.}%
\else
\language=\csname l@#1\endcsname
\fi
#2}}
\providecommand{\BIBdecl}{\relax}
\BIBdecl

\bibitem{JM-JB-JB:11}
J.~Machowski, J.~Bialek, and J.~Bumby, \emph{Power System Dynamics: Stability
  and Control}.\hskip 1em plus 0.5em minus 0.4em\relax John Wiley \& Sons,
  2011.

\bibitem{CZ-UT-LN-SL:14}
C.~Zhao, U.~Topcu, N.~Li, and S.~Low, ``Design and stability of load-side
  primary frequency control in power systems,'' \emph{IEEE Transactions on
  Automatic Control}, vol.~59, no.~5, pp. 1177--1189, 2014.

\bibitem{NL-CL-ZC-SHL:13}
N.~Li, L.~Chen, C.~Zhao, and S.~H. Low, ``Connecting automatic generation
  control and economic dispatch from an optimization view,'' in \emph{Proc.
  {A}merican {C}ontrol {C}onference}, Portland, OR, USA, 2014, pp. 735--740.

\bibitem{EM-CZ-SL:17}
E.~Mallada, C.~Zhao, and S.~Low, ``Optimal load-side control for frequency
  regulation in smart grids,'' \emph{IEEE Transactions on Automatic Control},
  vol.~62, no.~12, pp. 6294--6309, 2017.

\bibitem{EM-SL:14}
E.~Mallada and S.~H. Low, ``Distributed frequency-preserving optimal load
  control,'' \emph{IFAC Proceedings Volumes}, vol.~47, no.~3, pp. 5411--5418,
  2014.

\bibitem{SY-LC:14}
S.~You and L.~Chen, ``Reverse and forward engineering of frequency control in
  power networks,'' in \emph{Proc. {IEEE} Conference on Decision and Control},
  2014, pp. 191--198.

\bibitem{XZ-NL-AP:15}
X.~Zhang, N.~Li, and A.~Papachristodoulou, ``Achieving real-time economic
  dispatch in power networks via a saddle point design approach,'' in
  \emph{{IEEE} Power \& Energy Society General Meeting}, 2015.

\bibitem{TS-CDP-AVDS:17}
T.~Stegink, C.~De~Persis, and A.~van~der Schaft, ``A unifying energy-based
  approach to stability of power grids with market dynamics,'' \emph{{IEEE}
  Transactions on Automatic Control}, vol.~62, no.~6, pp. 2612--2622, 2017.

\bibitem{JWSP-FD-FB:13}
J.~W. Simpson-Porco, F.~D{\"o}rfler, and F.~Bullo, ``Synchronization and power
  sharing for droop-controlled inverters in islanded microgrids,''
  \emph{Automatica}, vol.~49, no.~9, pp. 2603--2611, 2013.

\bibitem{JWSP-QS-FD-JMV-JMG-FB:13s}
J.~W. Simpson-Porco, Q.~Shafiee, F.~D{\"o}rfler, J.~M. Vasquez, J.~M. Guerrero,
  and F.~Bullo, ``Secondary frequency and voltage control of islanded
  microgrids via distributed averaging,'' \emph{{IEEE} Transactions on
  Industrial Electronics}, vol.~62, no.~11, pp. 7025--7038, 2015.

\bibitem{FD-JWSP-FB:16}
F.~D{\"o}rfler, J.~W. Simpson-Porco, and F.~Bullo, ``Breaking the hierarchy:
  Distributed control and economic optimality in microgrids,'' \emph{{IEEE}
  Transactions on Control of Network Systems}, vol.~3, no.~3, pp. 241--253,
  2016.

\bibitem{FD-JWSP-FB:13y}
------, ``Plug-and-play control and optimization in microgrids,'' in
  \emph{Proc. {IEEE} Conference on Decision and Control}, Los Angeles, CA, USA,
  Dec. 2014, pp. 211--216.

\bibitem{QS-JMG-JCV:14}
Q.~Shafiee, J.~M. Guerrero, and J.~C. Vasquez, ``Distributed secondary control
  for islanded microgrids-a novel approach,'' \emph{{IEEE} Transactions on
  Power Electronics}, vol.~29, no.~2, pp. 1018--1031, 2014.

\bibitem{MA-DVD-KHJ-HS:13}
M.~Andreasson, D.~V. Dimarogonas, K.~H. Johansson, and H.~Sandberg,
  ``Distributed vs. centralized power systems frequency control,'' in
  \emph{Proc. European Control Conference}, 2013, pp. 3524--3529.

\bibitem{MA-DVD-HS-KHJ:14}
M.~Andreasson, D.~V. Dimarogonas, H.~Sandberg, and K.~H. Johansson,
  ``Distributed pi-control with applications to power systems frequency
  control,'' in \emph{Proc. American Control Conference}, 2014, pp. 3183--3188.

\bibitem{ST-MB-CDP:14}
S.~Trip, M.~B{\"u}rger, and C.~De~Persis, ``An internal model approach to
  frequency regulation in inverter-based microgrids with time-varying
  voltages,'' in \emph{Proc. {IEEE} Conference on Decision and Control}, 2014,
  pp. 223--228.

\bibitem{AB-FLL-AD:14}
A.~Bidram, F.~L. Lewis, and A.~Davoudi, ``Distributed control systems for
  small-scale power networks: Using multiagent cooperative control theory,''
  \emph{{IEEE} Control Systems Magazine}, vol.~34, no.~6, pp. 56--77, 2014.

\bibitem{CZ-EM-FD:15}
C.~Zhao, E.~Mallada, and F.~D{\"o}rfler, ``Distributed frequency control for
  stability and economic dispatch in power networks,'' in \emph{Proc. American
  Control Conference}, 2015, pp. 2359--2364.

\bibitem{FD-SG:16}
F.~D{\"o}rfler and S.~Grammatico, ``Amidst centralized and distributed
  frequency control in power systems,'' in \emph{Proc. American Control
  Conference}, 2016, pp. 5909--5914.

\bibitem{RM-SD-BBC:12}
R.~Mudumbai, S.~Dasgupta, and B.~B. Cho, ``Distributed control for optimal
  economic dispatch of a network of heterogeneous power generators,''
  \emph{{IEEE} Transactions on Power Systems}, vol.~27, no.~4, pp. 1750--1760,
  2012.

\bibitem{JWSP-BKP-NM-FD:18}
J.~W. Simpson-Porco, B.~K. Poolla, N.~Monshizadeh, and F.~Dorfler,
  ``Input-output performance of linear-quadratic saddle-point algorithms with
  application to distributed resource allocation problems,'' \emph{IEEE
  Transactions on Automatic Control}, {T}o appear. arXiv preprint
  arXiv:1803.02182.

\bibitem{JWSP-BKP-NM-FD:16}
J.~W. Simpson-Porco, B.~K. Poolla, N.~Monshizadeh, and F.~D{\"o}rfler,
  ``Quadratic performance of primal-dual methods with application to secondary
  frequency control of power systems,'' in \emph{Proc. {IEEE} Conference on
  Decision and Control}, 2016, pp. 1840--1845.

\bibitem{ARB:09}
A.~R. Bergen, \emph{Power Systems Analysis}.\hskip 1em plus 0.5em minus
  0.4em\relax Pearson Education India, 2009.

\bibitem{GED-FP:00}
G.~E. Dullerud and F.~Paganini, \emph{A Course in Robust Control Theory}, ser.
  Texts in Applied Mathematics.\hskip 1em plus 0.5em minus 0.4em\relax
  Springer, 2000, no.~36.

\bibitem{ET-BB-DFG:15}
E.~Tegling, B.~Bamieh, and D.~F. Gayme, ``The price of synchrony: Evaluating
  the resistive losses in synchronizing power networks,'' \emph{{IEEE}
  Transactions on Control of Network Systems}, vol.~2, no.~3, pp. 254--266,
  2015.

\bibitem{ET-MA-JWSP-HS:15d}
E.~Tegling, M.~Andreasson, J.~W. Simpson-Porco, and H.~Sandberg, ``Improving
  performance of droop-controlled microgrids through distributed
  {PI}-control,'' in \emph{Proc. {A}merrican {C}ontrol {C}onference}, Boston,
  MA, USA, Jul. 2016, pp. 2321--2327.

\bibitem{BKP-SB-FD:16}
B.~K. Poolla, S.~Bolognani, and F.~D{\"o}rfler, ``Placing rotational inertia in
  power grids,'' in \emph{Proc. American Control Conference}, 2016, pp.
  2314--2320.

\bibitem{ERAW-YJ-CZ-EM-CDP-FD:17}
E.~Weitenberg, Y.~Jiang, C.~Zhao, E.~Mallada, C.~De~Persis, and F.~D{\"o}rfler,
  ``Robust decentralized secondary frequency control in power systems: Merits
  and trade-offs,'' \emph{IEEE Transactions on Automatic Control}, November
  2017, {In press. DOI: 10.1109/TAC.2018.2884650. Available at
  {https://arxiv.org/abs/1711.07332}}.

\bibitem{CAD-ESK:69}
C.~A. Desoer and E.~S. Kuh, \emph{Basic Circuit Theory}.\hskip 1em plus 0.5em
  minus 0.4em\relax IEEE Press, 1969.

\bibitem{JS-DZ-RO-AMS-TS-JR:16}
J.~Schiffer, D.~Zonetti, R.~Ortega, A.~M. Stankovi{\'c}, T.~Sezi, and
  J.~Raisch, ``A survey on modeling of microgrids-from fundamental physics to
  phasors and voltage sources,'' \emph{Automatica}, vol.~74, pp. 135--150,
  2016.

\bibitem{VDB-JMH-AO-JNT:05}
V.~D. Blondel, J.~M. Hendrickx, A.~Olshevsky, and J.~N. Tsitsiklis,
  ``Convergence in multiagent coordination, consensus, and flocking,'' in
  \emph{Proc. {IEEE} Conference on Decision and Control}, 2005, pp. 2996--3000.

\end{thebibliography}

\section*{Appendix}

\begin{pfof}{Theorem~\ref{Thm:CAPIPerformance}}
We note that the state matrix is Hurwitz following \cite{FD-SG:16, JWSP-FD-FB:13} (consider the candidate potential function $S(\varphi, \omega, \mu)=\frac{1}{2}\varphi^\top \varphi +\frac{1}{2} m\,\omega^\top \omega + \frac{n\tau_\mu}{2k} \mu^2$). Further, under the given assumptions, for the states $(\varphi, \omega, \mu)$, the state matrix
\begin{align}
\label{Eq:CAPI-A}
A_{\rm agg} & =\begin{bmatrix}
\vzeros[] & E^\top & \vzeros[n-1]  \\[0.1cm]
-\frac{1}{m}E & -\frac{d}{m} I_n & -\frac{1}{m}K^{-1}\vones[n] \\[0.1cm]
\vzeros[n-1]^{\top} & \frac{1}{n\tau_\mu}\vones[n]^{\top}  & 0 \\
\end{bmatrix},
\end{align}
and the output matrix is ${C_{\rm agg}} = [\vzeros[n-1]^\top \, \vzeros[n]^\top\, -K^{-\frac{1}{2}}\vones[n]]$. 
One can verify by direct calculation that 
\begin{equation*}
X_{\rm agg}=\begin{bmatrix}
\vzeros[] & \vzeros[] & \vzeros[n-1]  \\[0.1cm]
\vzeros[] & z \vones[n] \vones[n]^{\top} & \frac{m}{2}\vones[n] \\[0.1cm]
\vzeros[n-1]^{\top} & \frac{m}{2}\vones[n]^{\top}  & \beta
\end{bmatrix},
\end{equation*}
is a positive semidefinite solution of the Lyapunov equation $X\, A_{\rm agg}+A_{\rm agg}^\top\, X+ C_{\rm agg}^\top \,C_{\rm agg}=\vzeros[]$, where $\beta\!=\! n\tau_{\mu}(d/2-z(\vones[n]^{\top} K^{-1} \vones[n])/m)$ and $z={m^2}/({2nd\tau_\mu})$. The $\mathcal{H}_2$ norm can be computed from \eqref{Eq:H2Obsv} as
\begin{align*}
\|G\|_{\mathcal{H}_2}^2 &= \mathrm{Tr}([\vzeros[]\,\,\tfrac{b}{m}I_n\,\,\vzeros[]]\, X_{\rm agg} [\vzeros[]\,\,\tfrac{b}{m}I_n\,\,\vzeros[]]^{\top})\\
&= z\frac{b^2}{m^2}\mathrm{Tr}(\vones[n]\vones[n]^{\top}) = z \frac{b^2}{m^2}n= \frac{b^2}{2 \tau_\mu}\frac{1}{d},
\end{align*}
which is independent of the system size.
\end{pfof}
\smallskip

\begin{pfof}{Theorem~\ref{Thm:PrimalDualFeedbackPerformance}}
We only provide a sketch of the proof due to space constraints. The approach developed in \cite{ET-BB-DFG:15,ET-MA-JWSP-HS:15d, JWSP-BKP-NM-FD:18} is followed here. Let $E = U\Sigma V^{\top}$ denote the singular value decomposition of the incidence matrix $E$, where $U \in\! \real^{n \times n}$ and $V \in\! \real^{|\mathcal{E}|\times|\mathcal{E}|}$ are orthogonal matrices and $\Sigma = \begin{bmatrix}\overline{\Sigma} & \vzeros[n \times (|\mathcal{E}|-n)]\end{bmatrix} \in \real^{n \times |\mathcal{E}|}$, where $\overline{\Sigma} = \mathrm{diag}(\sigma_{n}, \ldots, \sigma_2, 0)$ is positive semidefinite. We consider the coordinate transformation $(\hat{\varphi},\hat{\omega},\hat{\mu},\hat{\nu}) = (V^{\top}\varphi,U^{\top}\omega,U^{\top}\mu,V^{\top}\nu)$ along with the change of input $\hat{\eta} = U^{\top}\eta$ and output $\hat{y} = U^{\top}y$. Applying this transformation, we obtain the new system
\begin{equation}
\label{Eq:PrimalDual_Feedback_IO_Transformed}
\renewcommand*{\arraystretch}{1.2}
\begin{aligned}
\begin{bmatrix}
\dot{\hat{\varphi}}\\
\dot{\hat{\omega}}\\
\dot{\hat{\mu}}\\
\dot{\hat{\nu}}
\end{bmatrix}
&\!\!=\!\!
\underbrace{\begin{bmatrix}
\vzeros[] & \Sigma^\top & \vzeros[] & \vzeros[] \\ 
-\frac{1}{m}\Sigma & -\frac{d}{m}I_n & -\frac{1}{mk}I_n & \vzeros[] \\ 
\vzeros[] & \frac{\alpha}{\tau_\mu k}I_n & -\frac{1}{\tau_\mu k}I_n & -\frac{1}{\tau_\mu}\Sigma\\
\vzeros[] & \vzeros[] & \frac{1}{\tau_\nu}\Sigma^{\top} & \vzeros[]
\end{bmatrix}}_{A_\mathrm{pd}(\alpha)}\!\!\begin{bmatrix}\hat{\varphi} \\ \hat{\omega} \\ \hat{\mu} \\ \hat{\nu}\end{bmatrix}\!\!+\! \!
\underbrace{\begin{bmatrix}
\vzeros[]\\
\frac{b}{m}I_n\\
\frac{b}{\tau_{\mu}}I_n\\
\vzeros[]
\end{bmatrix}}_{B_\mathrm{pd}}\!\hat{\eta} \\
\hat{y} &= {-}K^{-\frac{1}{2}}\hat{\mu}, \, C_\mathrm{pd}=[\vzeros[]^\top\, \vzeros[]^\top\, -K^{-\frac{1}{2}}\, \vzeros[]^\top].
\end{aligned}
\end{equation}
It can be shown that the transformed system is internally stable (consider the candidate potential function $S(\hat\varphi, \hat\omega, \hat\mu, \hat\nu)=\frac{1}{2}\hat\varphi^\top \hat\varphi +\frac{1}{2}m\,\hat\omega^\top \hat\omega + \frac{\tau_\mu}{2\alpha} \hat\mu^\top \hat\mu +\frac{\tau_\nu}{2\alpha} \hat\nu^\top \hat\nu$). We recall that the $\mathcal{H}_2$ norm is invariant under unitary (orthogonal) transformations of inputs and outputs. Furthermore, for $X_\mathrm{pd}(\alpha)\!=\!\frac{1}{2}\,\mathrm{blkdiag}(\alpha I_{|\mathcal{E}|}, \alpha m I_n, \tau_{\mu} I_n, \tau_{\nu} I_{|\mathcal{E}|} )$,
the following inequality holds:
\[X_\mathrm{pd}(\alpha)\,A_\mathrm{pd}(\alpha)+A_\mathrm{pd}(\alpha)^\top X_\mathrm{pd}(\alpha)+Q^i_\mathrm{pd}\preceq\vzeros[],\]
where $Q_\mathrm{pd}= {C_{\rm pd}}^\top {C_{\rm pd}}$.
Therefore, by \cite[Chapter 4.7]{GED-FP:00}, $X_\mathrm{pd}(\alpha)$ is the generalized Gramian and we have that the squared $\mathcal{H}_2$ norm of the overall system is upper-bounded as
\[\|G(\alpha)\|_{\mathcal{H}_2}^2 \!\leq \!\text{Tr}(B_{\mathrm{pd}}^\top X_\mathrm{pd}(\alpha) B_\mathrm{pd})\!=\!\|G(0)\|_{\mathcal{H}_2}^2+\frac{b^2\alpha n}{2m}.\]\
We note here that the upper-bound is linear in the gain $\alpha$, and we recover equality for $\alpha=0$. Furthermore, in the limit of large generators, i.e., $m\to \infty$ or small networks $n \to 0$, we arrive back at the open-loop performance.
\end{pfof}

\smallskip

\begin{pfof}{Theorem~\ref{Thm:DAPI_Performance}}
We note that as $L \in \real^{n \times n}$ is symmetric, there exists an orthogonal matrix $O \in \real^{n \times n}$ such that $L = O\Lambda O^{\top}$, where $\Lambda = \mathrm{diag}(0,\lambda_2,\ldots,\lambda_n)$ is the positive semidefinite diagonal matrix of eigenvalues of $L$. Consider now the change of coordinates $(\hat{\theta},\hat{\omega},\hat{p}) = (O{\theta},O{\omega},O{p})$ along with the change of input $\hat{\eta} = O\eta$ and output $\hat{y} = Oy$. In these new coordinates, all blocks in the system matrices are diagonalized, and the system becomes the $n$ decoupled systems. If the input-output mapping of the transformed $i$th subsystem is denoted by $\hat{G}_i$, then it is straightforward to verify that $\|{G}\|_{\mathcal{H}_2}^2 = \sum_{i=1}^n \|\hat{G}_i\|_{\mathcal{H}_2}^2$, so our problem reduces to computing the $\mathcal{H}_2$ norm of each decoupled subsystem $\hat{G}_i$. 

For $i \in \{2,\ldots,n\}$, one may verify without much difficulty that the $i$th decoupled subsystem is internally stable and hence subsystem $\hat{G}_i$ of the transformed system is input-output stable with finite norm. For $\lambda_1=0$, the marginally stable mode is unobservable from the output, and hence also input-output stable with a finite norm. For each subsystem, the Lyapunov equation \eqref{Eq:ObsvGram} is a $3\times 3$ matrix equality which can be parameterized, expanded, and solved by hand. The results of Theorem \ref{Thm:DAPI_Performance} then immediately follow by applying \eqref{Eq:H2Obsv} to each subsystem, summing the results, and taking the appropriate limits for (i) and (ii).
\end{pfof}

\smallskip

\end{document}